\documentclass[12pt,a4paper]{amsart}
\usepackage{a4wide}
\usepackage[english]{babel}
\usepackage[T2A]{fontenc}
\usepackage[utf8]{inputenc} 
\usepackage{amsfonts}
\usepackage{amssymb, amsthm, amscd}
\usepackage{amsmath}
\usepackage{mathtools}
\usepackage{needspace}
\usepackage{etoolbox}
\usepackage{lipsum}
\usepackage{comment}
\usepackage{cmap}
\usepackage[pdftex]{graphicx}
\usepackage[unicode]{hyperref}
\usepackage[matrix,arrow,curve]{xy}
\usepackage[usenames,dvipsnames]{xcolor}
\usepackage{colortbl}
\usepackage{textcomp}
\usepackage{cite}
\usepackage{euscript}
\usepackage[numbers]{natbib}

\newcommand{\Z}{\mathbb{Z}}

\renewcommand{\S}{\mathbb{S}}

\newcommand{\Pm}{\mathbb{P}}
\renewcommand{\P}{\mathfrak{p}}
\renewcommand{\frq}{\mathfrak{q}}

\newcommand{\A}{\mathfrak{A}}

\renewcommand{\O}{\EuScript{O}}

\newcommand{\ovl}{\overline}

\newcommand{\eps}{\varepsilon}

\newcommand{\sub}{\subseteq}

\newcommand{\Ker}{\mathop{\mathrm{Ker}}\nolimits}

\newcommand{\rk}{\mathop{\mathrm{rk}}\nolimits}
\newcommand{\Cl}{\mathop{\mathrm{Cl}}\nolimits}

\renewcommand{\ge}{\geqslant}
\renewcommand{\le}{\leqslant}
\newcommand{\sm}{\setminus}

\newcommand{\SL}{\mathop{\mathrm{SL}}\nolimits}

\newcommand{\Gal}{\mathop{\mathrm{Gal}}\nolimits}

\newcommand{\tc}{\text{,}}
\newcommand{\tp}{\text{.}}
\renewcommand{\tilde}{\widetilde}

\renewcommand{\b}{\mathfrak{b}}

\newcommand{\rsymb}[2]{\genfrac{(}{)}{}{}{#1}{#2}}

\DeclareMathOperator{\res}{res}
\DeclareMathOperator{\ord}{ord}
\DeclareMathOperator{\bad}{bad}
\DeclareMathOperator{\lcm}{lcm}

\theoremstyle{plain}
\newtheorem{thm}{Theorem}[section]
\newtheorem{lem}[thm]{Lemma}
\newtheorem{prop}[thm]{Proposition}

\theoremstyle{definition}

\newtheorem{cor}[thm]{Corollary}

\theoremstyle{remark}
\newtheorem{rem}[thm]{Remark}

\AtBeginEnvironment{thm}{\begin{samepage}}
	\AtEndEnvironment{thm}{\end{samepage}}
\AtBeginEnvironment{lem}{\begin{samepage}}
	\AtEndEnvironment{lem}{\end{samepage}}
\AtBeginEnvironment{st}{\begin{samepage}}
	\AtEndEnvironment{st}{\end{samepage}}
\AtBeginEnvironment{crit}{\begin{samepage}}
	\AtEndEnvironment{crit}{\end{samepage}}
\AtBeginEnvironment{ax}{\begin{samepage}}
	\AtEndEnvironment{ax}{\end{samepage}}
\AtBeginEnvironment{defn}{\begin{samepage}}
	\AtEndEnvironment{defn}{\end{samepage}}
\AtBeginEnvironment{cor}{\begin{samepage}}
	\AtEndEnvironment{cor}{\end{samepage}}
\AtBeginEnvironment{note}{\begin{samepage}}
	\AtEndEnvironment{note}{\end{samepage}}
\AtBeginEnvironment{prop}{\begin{samepage}}
	\AtEndEnvironment{prop}{\end{samepage}}

\makeatletter
\def\@settitle{\begin{center}%
		\baselineskip14\p@\relax
		\bfseries
		\@title
	\end{center}%
}

\def\@evenhead{\hfil\sc Pavel Gvozdevsky\hfil}
\def\@oddhead{\hfil\sc Width\hfil}
\makeatother

\title{Width of $\SL(n,\O_S,I)$}


\author{Pavel Gvozdevsky}
\date{}
\address{Chebyshev Laboratory, St. Petersburg State University, 14th Line V.O., 29B, Saint Petersburg 199178 Russia}
\email{gvozdevskiy96@gmail.com}
\thanks{Research is supported by Russian Science Foundation grant (project \textnumero 22-21-00257).}

\begin{document}

\maketitle

\begin{abstract}
	We give an estimate for the width of the congruence subgroup $\SL(n,\O_S,I)$ in Tits--Vaserstein generators, where $\O_S$ is a localisation of the ring of integers in a number field $K$. We assume that either $K$ has a real embedding, or the ideal $I$ is prime to the number of roots of unity in $K$.
\end{abstract}

\section{Introduction}

Given a group $G$ with generating set $X$, the width of $G$ in generators from $X$ is a minimal number $N$ such that any element of $G$ is a product of at most $N$ elements from $X$. When the width is finite, we say that $G$ admits bounded generation with respect to $X$.

For special linear group, or more generally for Chevalley groups bounded generation with respect to the elementary generators is known for certain classes of rings. For example this holds for Dedekind domains of arithmetic type, see \cite{CarterKeller}, \cite{CarterKellerZ}, \cite{Morris}, \cite{TavgenChevalley}, \cite{TavgenTwisted},\cite{MorganRapinchukSury},\cite{ErovenkoRapinchuk2001},\cite{ErovenkoRapinchuk2006},\cite{KunPlotVavBounded}, \cite{TrostFiniteIndex}, \cite{TrostChevalleyStrong}, \cite{TrostChevalley},\cite{TrostSL2},\cite{TrostSp4}, \cite{TrostSympl}. These results are of great value, for example they are connected to the {\it congruence subgroup property}, see \cite{Lubotzky},\cite{PlatonovRapinchuk}; to Margulis--Zimmer conjecture, see \cite{ShalomWillis}; and have applications in logic, see, for example, \cite{AvniLubMeiri}, \cite{AvniMeiri}.

In \cite{SinchukSmolensky}, it was proven that the principal congruence subgroup $G(\Phi,\O_S,I)$ of a classical Chevalley group with $\rk\Phi\le 2$ over a Dedekind domain $\O_S$ of arithmetic type has finite width in Tits--Vaserstein generators, provided the fraction field of $\O_S$ has a real embedding. Also this result can be deduced from \cite{TrostFiniteIndex}. However, either way the proof relies on results from \cite{TavgenThesis}, which on its turn are not constructive and do not allow to obtain an explicit estimate of the width in question. In Section \ref{finalremarks}, we show how one can one can obtain such a proof directly from the result of \cite{CarterKeller}.

In the main part of the present paper we prove an effective version of the result from \cite{SinchukSmolensky} for special linear group, i.e. the width of $\SL(n,\O_S,I)$ will be estimated explicitly.

First we consider the double $\tilde{\O}$ of the ring $\O$ with respect to the ideal $I$, where $\O$ is the ring of integers in a number field, and adopt the technique from \cite{CarterKeller} in order to give an estimate for how many elementary transvections from $\SL(3,\tilde{\O})$ is enough to present any matrix from $\SL(2,\tilde{\O})$ (Theorem~\ref{main}). Then we generalise this result to a localisation $\O_S$ of the ring $\O$ (Corollary~\ref{loc}). After that we use these results to estimate the width of $\SL(n,\O_S,I)$ in Tits--Vaserstein generators (Corollary~\ref{WidthTitsVas}). Finally, considering two ideals $A$,$B\unlhd \O_S$, we use the proof of Lemma~7 in \cite{VavCuriouser} in order to estimate the width of $\SL(n,\O_S,AB)$ in generators $[t_{i,j}(a),t_{h,k}(b)]$, where $1\le i\ne j\le n$, $1\le h\ne k\le n$, $a\in A$ and $b\in B$.

The paper is organised af follows. In Section~\ref{PreliminariesAndNotation} we give the necessary preliminaries and introduce basic notation. In Section~\ref{statement} we formulate Theorem~\ref{main} and prove its corollaries as described above. In Section~\ref{proof} we give the proof of Theorem~\ref{main}. In Section~\ref{finalremarks} we give some final remarks.

I am grateful to Nikolai Vavilov for optimising the proof of Corollary~\ref{TitsVaserstein}.

\section{Preliminaries and notation}
\label{PreliminariesAndNotation}
\subsection{Special linear group}

For a commutative ring $R$ we denote by $\SL(n,R)$ the group of $n\times n$ matrices over $R$ with determinant one. By $t_{i,j}(\xi)$, where $1\le i\ne j\le n$ and $\xi\in R$, we denote an elementary transvection, i.e. the matrix that has $1$ in all the diagonal positions, $\xi$ in the position $(i,j)$ and $0$ in all the remaining positions. By $E(n,R)$ we denote the elementary group, i.e. the subgroup of $\SL(n,R)$ generated by all the elementary transvections. 

\subsection{Congruence subgroups}

For an ideal $I\unlhd R$, we denote by $\SL(n,R,I)$ the principal congruence subgroup, i.e the kernel of the reduction homomorphism $\SL(n,R)\to \SL(n,R/I)$ induced by the projection $R\to R/I$. By $E(n,I)$ we denote the group generated by all the elementary transvections $t_{i,j}(\xi)$ with $\xi\in I$; and by $E(n,R,I)$ we denote the elementary congruence subgroup, i.e. the normal closure of $E(n,I)$ in $E(n,R)$.

As a group $E(n,R,I)$ is generated by Tits--Vaserstein generators $$\{t_{i,j}(\xi)^{t_{j,i}(\zeta)}\colon \xi\in I\tc\; \zeta\in\O_S\}\tc$$ see, for example, Theorem 2 in \cite{VasersteinChevalley}. 

In the present paper, we consider the localisation $\O_S$ of the ring $\O$ of integers in a number field $K$. We assume that either $K$ has a real embedding, or the ideal $I$ is prime to $m$, where $m$ be the number of roots of unity in $K$. Under these assumptions it follows from \cite{BassMilnorSerre} that $E(n,\O_S,I)=\SL(n,\O_S,I)$. Corollary~\ref{WidthTitsVas} of the present paper gives an estimate for the width of $\SL(n,\O_S,I)$ in Tits--Vaserstein generators.

\subsection{Double of a ring}

Let $R$ be a commutative ring, and $I\unlhd R$ be an ideal. The {\it double} of the ring $R$ with respect to the ideal $I$ is a subring of $R\times R$ that consists of pairs with elements congruent to each other modulo $I$:

$$
\tilde{R}=\{(a,b)\in R\times R\colon a\equiv b\mod I\}\tp
$$

We will use the following lemma.

{\lem \label{quotient} Let $b=(b',b'')\in \tilde{R}$ be such that $b'$ {\rm (}and hence $b''${\rm )} is prime to $I$. Then the natural map
$$
\tilde{R}/b\tilde{R}\to R/b'R\times R/b''R
$$
is an isomorphism.}

\begin{proof}
	It is easy to see that $\Ker (\tilde{R}\to R/b'R)=b\tilde{R} + (0,I)$ and $\Ker (\tilde{R}\to R/b''R)=b\tilde{R} + (I,0)$. These two ideals are comaximal, because $\tilde{R}/(I,I)=R/I$ and by assumption the image of $b$ in $R/I$ is invertible. Hence we have
	$$
	b\tilde{R}\le (b\tilde{R} + (0,I))\cap (b\tilde{R} + (I,0))=(b\tilde{R} + (0,I)) (b\tilde{R} + (I,0))\le b\tilde{R}\tp
	$$
	So we have 
	$$
	b\tilde{R} + (0,I))\cap (b\tilde{R} + (I,0))=b\tilde{R}\tp
	$$
	The statement of the lemma now follows from the Chinese remainder theorem.
	\end{proof}

\section{Statement of the main theorem and proofs of corollaries.}
\label{statement}

Let $\O$ be the ring of integers in an algebraic number field $K$. Let $D$ be the discriminant of $K$ and $\Cl(K)$ be its class group. Let $m$ be the number of roots of unity in $K$. For any rational prime $p$ set $e_p=\ord_p(m)$, i.e. $m=\prod_{\{p\colon e_p>0\}} p^{e_p}$. Further for any rational prime $p$ we denote by $L_p$ the extension of $K$ obtained by adjoining a primitive $p^{e_p+1}$-th root of unity. Now set
$$
\S_{\bad}=\{p\in \Pm\colon p\mid D\text{ and } \gcd ([L_p\colon K], |\Cl(K)|)>1 \}\tc
$$
where $\Pm$ denotes the set of rational primes. Finally, set 
$$
\Delta=\max_{\delta_1+\delta_2+\delta_3=|\S_{\bad}|}\left(\sum_{i=1}^3\max(1,[\ln(\delta_i+1)/\ln 2])\right)\tc
$$
where maximum is taken over all triples of nonnegative integers $\delta_1$,$\delta_2$,$\delta_3$ with $\delta_1+\delta_2+\delta_3=|\S_{\bad}|$.

\medskip

The main result of the present paper is the following theorem.

{\thm\label{main} In the notation above, let $I$ be a non zero ideal in $\O$. Suppose that either $K$ has a real embedding, or $I$ is prime to $m$. Let $\tilde{\O}$ be the double of the ring $\O$ with respect to the ideal $I$. Then for any matrix $\left(\begin{smallmatrix} a& b\\ c& d \end{smallmatrix}\right)\in \SL(2, \tilde{\O})$ the matrix
$$
\begin{pmatrix}
	a & b & 0\\
	c & d & 0\\
	0 & 0 & 1
\end{pmatrix}\in \SL(3,\tilde{\O})
$$
is a product of at most $68\Delta+4$ elementary transvections.}

\medskip

We prove this theorem in Section \ref{proof}. Now we deduce the corollary for localisations of the ring $\O$. 

{\cor\label{loc} In the notation above, let $S$ be a multiplicative system in $\O$, let $\O_S=\O[S^{-1}]$, and $I_S$ be a non zero ideal in $\O_S$. Suppose that either $K$ has a real embedding, or $I_S$ is prime to $m$. Let $\tilde{\O_S}$ be the double of the ring $\O_S$ with respect to the ideal $I_S$. Then for any matrix $\left(\begin{smallmatrix} a& b\\ c& d \end{smallmatrix}\right)\in \SL(2, \tilde{\O_S})$ the matrix
	$$
	\begin{pmatrix}
		a & b & 0\\
		c & d & 0\\
		0 & 0 & 1
	\end{pmatrix}\in \SL(3,\tilde{\O_S})
	$$
	is a product of at most $68\Delta+8$ elementary transvections.} 

\begin{proof}
	Let $I$ be the preimage of $I_S$ in $\O$. Note that any element $s\in S$ is prime to $I$. Indeed, let $I=\prod \P_i^{k_i}$ be the decomposition of $I$ into a product of primes, and assume that $s$ belongs to one of the $\P_i$, say $\P_1$. Then the localisation homomorphism maps $\prod_{i\ne 1} \P_i^{k_i}$ to $I_S$; hence $\prod_{i\ne 1} \P_i^{k_i}\le I$, which is a contradiction.
	
	It follows that, if $I_S$ is prime to $m$ in $\O_S$, then $I$ is prime to $m$ in $\O$. Indeed, in this case the ideal $m\O+I$ must contain an element $s\in S$; and since $m\O+I$ also contains $I$, it follows that $m\O+I=\O$.
	
	Therefore, Theorem \ref{main} can be applied to the ring $\O$ and the ideal $I$. It is easy to see that $\tilde{\O_S}=\tilde{\O}[S^{-1}]$, where $S$ is embeded into $\tilde{\O}$ diagonally; and $\tilde{\O}$ maps to $\tilde{\O_S}$ injectively. So it remains to prove that any matrix $\left(\begin{smallmatrix} a_1& b_1\\ c_1& d_1 \end{smallmatrix}\right)\in \SL(2, \tilde{\O_S})$ can be transformed to a matrix from $\SL(2,\tilde{\O})$ by 4 elementary transformations.
	
	We perform these transformations as follows. At the first step, we multiply our matrix from the right by a suitable transvection, so that in the new matrix  $\left(\begin{smallmatrix} a_2& b_1\\ * & * \end{smallmatrix}\right)$ the entry $a_2$ became prime to $(I_S, I_S)$. This is possible because the ring $\tilde{\O}_S/(I_S,I_S)\simeq \O/I$ is semilocal. 
	
	Now let $a_2=(a_2',a_2'')$ and $b_1=(b_1',b_1'')$. Both $a_2'$ and $a_2''$ are prime to $I_S$, so in particular, they are non zero; hence $\O\to \O_s/a_2'\O_S$ is surjective. Hence we can choose $b_2'\in \O$ such that $b_2'\equiv b_1'\mod a_2'\O_S$. Moreover, since $a_2'$ is prime to $I_S$ and every element of $S$ is prime to $I$ in $\O$, it follows that the numerator of $a_2'$ is prime to $I$ in $\O$; hence changing $b_2'$ by a multiple of that numerator, we may assume that $b_2'\equiv 1\mod I$. Similarly, we choose $b_2''\in \O$ so that $b_2''\equiv b_1''\mod a_2''\O_S$ and $b_2''\equiv 1\mod I$. By construction, we have $b_2=(b_2',b_2'')\in \tilde{\O}$. Since $b_2'\equiv b_1'\mod a_2'\O_S$, $b_2''\equiv b_1''\mod a_2''\O_S$, and $a_2'$ is prime to $I_S$, it follows by Lemma~\ref{quotient} that $b_2\equiv b_1\mod a_2\tilde{\O_S}$. Therefore, at the second step, we can multiply our matrix from the right by a suitable transvection to transform it into the matrix $\left(\begin{smallmatrix} a_2& b_2\\ * & * \end{smallmatrix}\right)$.
	
	Now we factor $b_2'\O=\b_1'\b_2'$, such that no prime divisors of $\b_1'$ meet $S$ and such that all prime divisors of $\b_2'$ meet $S$. Similarly, we factor $b_2''\O=\b_1''\b_2''$. Further choose, $v',v''\in \O$ such that $v'\equiv a_2'\mod b_2'\O_S$ and $v''\equiv a_2''\mod b_2''\O_S$. By the Chinese remainder theorem we may now choose $a_3',a_3''\in\O$ such that 
	\begin{align*}
		&a_3'\equiv v'\mod \b_1'\O\tc\qquad a_3'\equiv 1\mod \b_2'I\tc\\
		&a_3''\equiv v''\mod \b_1''\O\tc\qquad a_3''\equiv 1\mod \b_2''I\tp
	\end{align*}
By construction, we have $a_3=(a_3',a_3'')\in\tilde{\O}$. Since $a_3'\equiv a_2'\mod b_2'\O_S$, $a_3''\equiv a_2''\mod b_2''\O_S$, and $b_2'$ is prime to $I_S$, it follows by Lemma~\ref{quotient} that $a_3\equiv a_2\mod b_2\tilde{\O_S}$. Therefore, at the third step, we can multiply our matrix from the right by a suitable transvection to transform it into the matrix $\left(\begin{smallmatrix} a_3& b_2\\ * & * \end{smallmatrix}\right)$.

In addition, by construction, we have $a_3'\O+b_2'\O=a_3''\O+b_2''\O=\O$. Thus $a_3'$ is invertible modulo $b_2'\O$ and  $a_3''$ is invertible modulo $b_2''\O$. It follows then by Lemma~\ref{quotient} that $a_3$ is invertible modulo $b_2\tilde{\O}$, i.e. $a_3\tilde{\O}+b_2\tilde{\O}=\tilde{\O}$. Therefore, at the forth step, we can multiply our matrix from the left by a suitable transvection to transform it into the matrix from $\SL(2,\tilde{\O})$.
\end{proof}

Now we show that the results above allow to estimate the width of the congruence subgroup $\SL(n,\O_S,I)$. First we consider matrices from the image of the embedding $\SL(2,\O_S,I)\to \SL(3,\O_S,I)$.

{\cor\label{conjugates} In the notation above, let $S$ be a multiplicative system in $\O$, let $\O_S=\O[S^{-1}]$, and $I$ be a non zero ideal in $\O_S$. Suppose that either $K$ has a real embedding, or $I$ is prime to $m$. Then for any matrix $\left(\begin{smallmatrix} a& b\\ c& d \end{smallmatrix}\right)\in \SL(2,\O_S, I)$ the matrix
	$$
	g=\begin{pmatrix}
		a & b & 0\\
		c & d & 0\\
		0 & 0 & 1
	\end{pmatrix}\in \SL(3,\O_S,I)
	$$
	is a product of at most $68\Delta+8$ elements of type $t_{i,j}(\xi)^h$, where $\xi\in I$ and $h\in \SL(3,\O_S)$.}

\begin{proof}
Set $N=68\Delta+8$. Consider the matrix
$$
\begin{pmatrix}
	(a,1) & (b,0) & (0,0)\\
	(c,0) & (d,1) & (0,0)\\
	(0,0) & (0,0) & (1,1)
\end{pmatrix}\in \SL(3,\tilde{\O_S})\tp
$$
By Corollary \ref{loc}, this matrix is a product of $N$ elementary transvections. Thus there exist elements $\zeta_1'$,$\ldots$,$\zeta_N'$,$\zeta_1''$,$\ldots$,$\zeta_N''\in\O_S$ such that 

$\bullet$ $\zeta_k'\equiv\zeta_k''\mod I$ for all $1\le k\le N$,

$\bullet$ $\prod_{k=1}^N t_{i_k,j_k}(\zeta_k')=g$,

$\bullet$ $\prod_{k=1}^N t_{i_k,j_k}(\zeta_k'')=e$,

for some indices $i_k$, $j_k$.

Set $\xi_k=\zeta_k'-\zeta_k''\in I$. Then we have

$$
g=\prod_{k=1}^N t_{i_k,j_k}(\zeta_k'')t_{i_k,j_k}(\xi_k)=\left(\prod_{k=1}^N t_{i_k,j_k}(\zeta_k'')\right) \left(\prod_{k=1}^N t_{i_k,j_k}(\xi_k'')^{h_k} \right)=\prod_{k=1}^N t_{i_k,j_k}(\xi_k)^{h_k}\tc
$$

where $h_k=\prod_{l=k+2}^N t_{i_l,j_l}(\zeta_l'')$.
\end{proof}

\medskip

{\cor\label{TitsVaserstein} Under the conditions of Corollary \ref{conjugates}, the matrix $g$ is a product of at most $24\cdot(68\Delta+8)=1632\Delta+192$ Tits--Vaserstein generators, i.e. elements of type $t_{i,j}(\xi)^{t_{j,i}(\zeta)}$, where $\xi\in I$, $\zeta\in\O_S$.}

\begin{proof}
	Due to Corollary \ref{conjugates}, it is enough to prove that any element of type $t_{i,j}(\xi)^h$ is a product of at most 24 elements of type $t_{i,j}(\xi)^{t_{j,i}(\zeta)}$.
	
	So let $g=t_{i,j}(\xi)^h$. Without loss of generality, we assume that $i=1$, $j=2$. 
	
	We use the technique from \cite{StepVavDecomp}. Let
	$$
	h=\begin{pmatrix}
				h_{1,1} & h_{1,2} & h_{1,3}\\
				h_{2,1} & h_{2,2} & h_{2,3}\\
				h_{3,1} & h_{3,2} & h_{3,3}
			\end{pmatrix}\qquad
		h^{-1}=\begin{pmatrix}
			h'_{1,1} & h'_{1,2} & h'_{1,3}\\
			h'_{2,1} & h'_{2,2} & h'_{2,3}\\
			h'_{3,1} & h'_{3,2} & h'_{3,3}
		\end{pmatrix}\tp
		$$
		Then we have
		$$
		t_{1,2}(\xi)=\prod_{k=1}^3 t_{1,2}(h'_{k,3}\xi h_{3,k})t_{1,3}(-h'_{k,3}\xi h_{2,k})\tp
		$$
		Thus it is enough to express $(t_{1,2}(h'_{k,3}\xi h_{3,k})t_{1,3}(-h'_{k,3}\xi h_{2,k}))^h$ as a product of 8 Tits--Vaserstein generators.  It follows from \cite{StepVavDecomp} that
		$$
		(t_{1,2}(h'_{k,3}\xi h_{3,k})t_{1,3}(-h'_{k,3}\xi h_{2,k}))^h=[t_{i,k}(\zeta_1)t_{j,k}(\zeta_2),t_{k,i}(\xi_1)t_{k,j}(\xi_2)]t_{k,i}(\xi_3)t_{k,j}(\xi_4)\tc
		$$
		where $\{i,j,k\}=\{1,2,3\}$, $\xi_1=h'_{k,3}\xi h'_{j,1}$, $\xi_2=-h'_{k,3}\xi h'_{i,1}$, $\xi_4=h'_{k,1}\xi_1$, $\xi_5=h'_{k,1}\xi_2$, $\zeta_1=h'_{i,1}$, $\zeta_2=h'_{j,1}$.
		
		Now we have
		$$
		(t_{1,2}(h'_{k,3}\xi h_{3,k})t_{1,3}(-h'_{k,3}\xi h_{2,k}))^h= (t_{k,i}(\xi_1)t_{k,j}(\xi_2))^{t_{i,k}(-\zeta_1)t_{j,k}(-\zeta_2)}t_{k,i}(\xi_3-\xi_1)t_{k,j}(\xi_4-\xi_2)\tp
		$$
		Thus it remains to express $(t_{k,i}(\xi_1)t_{k,j}(\xi_2))^{t_{i,k}(-\zeta_1)t_{j,k}(-\zeta_2)}$ as a product of 6 Tits--Vaserstein generators.
		
		We have
		$$
		(t_{k,i}(\xi_1)t_{k,j}(\xi_2))^{t_{i,k}(-\zeta_1)t_{j,k}(-\zeta_2)}=t_{k,i}(\xi_1)^{t_{i,k}(-\zeta_1)t_{j,k}(-\zeta_2)}t_{k,j}(\xi_2)^{t_{i,k}(-\zeta_1)t_{j,k}(-\zeta_2)}\tc
		$$
		where
		$$
		t_{k,i}(\xi_1)^{t_{i,k}(-\zeta_1)t_{j,k}(-\zeta_2)}=(t_{j,i}(\xi_1\zeta_2)t_{k,i}(\xi_1))^{t_{i,k}(-\zeta_1)}=t_{j,i}(\xi_1\zeta_2)t_{j,k}(-\xi_1\zeta_1\zeta_2)(t_{k,i}(\xi_1))^{t_{i,k}(-\zeta_1)}\tc
		$$
		and similarly
		$$
		t_{k,j}(\xi_2)^{t_{i,k}(-\zeta_1)t_{j,k}(-\zeta_2)}=t_{i,j}(\xi_2\zeta_1)t_{i,k}(-\xi_2\zeta_1\zeta_2)(t_{k,j}(\xi_2))^{t_{j,k}(-\zeta_2)}\tp
		$$
		That finishes the proof.
	\end{proof}

\medskip

Now we estimate the width of $\SL(n,\O_S,I)$ for arbitrary $n\ge 3$.

{\cor\label{WidthTitsVas} In the notation above, let $S$ be a multiplicative system in $\O$, let $\O_S=\O[S^{-1}]$, and $I$ be a non zero ideal in $\O_S$. Suppose that either $K$ has a real embedding, or $I$ is prime to $m$. Let $n\ge 3$. The the width of $\SL(n,\O_S,I)$ in Tits--Vaserstein generators $\{t_{i,j}(\xi)^{t_{j,i}(\zeta)}\colon \xi\in I\tc\; \zeta\in\O_S\}$ is at most $3n(n-1)/2+2n+1632\Delta+185$.}

\begin{proof}
	By Corollary 4.8 in \cite{SinchukSmolensky}, any every element of $\SL(n,\O_S,I)$ can be decomposed into a product of one element of $\SL(2,\O_S,I)$ and at most $3n(n-1)/2+2n-7$ Tits--Vaserstein generators. It remains to apply Corollary~\ref{TitsVaserstein} to that one element of $\SL(2,\O_S,I)$.
	\end{proof}

\smallskip

{\rem Similarly one can give an estimate for the width of a congruence subgroup of the split spin group. Unfortunately, at the time of this writing there are no proven analog of Corollary 4.8 in \cite{SinchukSmolensky} for Chevalley groups of type E; however, one can still give an estimate for the width of a corresponding congruence subgroup as follows. Since the double $\tilde{\O_S}$ of the ring $\O_S$ with respect to the ideal $I$ has Krull dimension equal to one, the proofs of the surjective $K_1$-stability as presented in \cite{SteinStability},\cite{PlotkinStability},\cite{PlotkinE7},\cite{GvozK1Surj} together with Corollary~\ref{loc} allow to estimate the width of the Chevalley group over $\tilde{\O_S}$ in elementary generators. This gives an estimate for the width of a congruence subgroup in conjugates of elementary generators. After that it is theoretically possible to estimate how many Tits-Vaserstein generators one need to express such a conjugate in a manner similar to the proof of Corollary~\ref{TitsVaserstein}. }

{\cor\label{WidthCommutators} In the notation above, let $S$ be a multiplicative system in $\O$, let $\O_S=\O[S^{-1}]$, let $A$ and $B$ be non zero ideals in $\O_S$. Suppose that either $K$ has a real embedding, or $AB$ is prime to $m$. Let $n\ge 3$. Then the group $\SL(n,\O_S,AB)$ is generated by commutators $[t_{i,j}(a),t_{h,k}(b)]$, where $1\le i\ne j\le n$, $1\le h\ne k\le n$, $a\in A$ and $b\in B$; and the width of $\SL(n,\O_S,AB)$ in these generators is at most  $12\cdot(3n(n-1)/2+2n+1632\Delta+185)=18n^2+6n+19584\Delta+2220$.}

\begin{proof}
	Due to Corollary~\ref{WidthTitsVas} it is enough to prove that any Tits--Vaserstein generator for $\SL(n,\O_S,AB)$ is a product of at most 12 such commutators; and that follows from the proof of Lemma~7 in \cite{VavCuriouser}.
\end{proof}

\smallskip

{\rem Similarly, if one calculate an estimate for the width of congruence subgroup in the spin group or the Chevalley group of type E in Tits--Vaserstein generators, then multiplying this  number by 12 one gets an estimate of the width in elementary commutators similar to those above.}

\section{The proof of the main theorem}
\label{proof}

For a commutative  ring $R$ and an ideal $J\unlhd R$ we denote by $\eps_R(J)$ the exponent of the multiplicative group $(R/J)^*$ with the convention $\eps_R(J)=0$ if this exponent is infinite. In other words, $\eps_R(J)$ is the smallest positive integer such that for any element $a\in R$ prime to $J$ we have $a^{\eps_R(J)}\equiv 1\mod J$, or zero if such a positive integer does not exists. For any $b\in R$ we set $\eps_R(b)=\eps_R(bR)$. Also we will omit the ring in the index if it is clear from the context. 

\medskip

The proof of the next proposition is the same as the proof of the main theorem in \cite{CarterKeller}, so we omit it here.

{\prop\label{general} Let $R$ be a commutative ring; let $k$, $\Delta$ and $m$ be positive integers; and let $\Omega\sub R^2$. Suppose that the two following conditions hold true.
\begin{enumerate}
\item Any row $(a_1,b_1)$, where $a_1,b_1\in R$ and $a_1R+b_1R=R$, can be transformed to a row $(a^m,b)$, where $(a,b)\in\Omega$, by $k$ elementary transvections from $\SL(2,R)$.

\item For any $(a,b)\in \Omega$ there exist $c_1,\ldots,c_{\Delta_0}\in R$ and $\gamma_1,\ldots,\gamma_{\Delta_0}\in \Z_{\ge 0}$, where $\Delta_0\le\Delta$, such that $\gcd(\gamma_1,\ldots,\gamma_{\Delta_0})=1$, and for any $1\le i\le {\Delta_0}$ we have $c_ib\equiv -1\mod aR$ and $\gcd(\eps(b),\eps(c_i))=m\gamma_i$.

\smallskip

Then for any matrix $\left(\begin{smallmatrix} a& b\\ c& d \end{smallmatrix}\right)\in \SL(2, R)$ the matrix
$$
\begin{pmatrix}
	a & b & 0\\
	c & d & 0\\
	0 & 0 & 1
\end{pmatrix}\in \SL(3,R)
$$
is a product of at most $68\Delta+k$ elementary transvections.
\end{enumerate} }

In order to proof Theorem \ref{main}, we apply Proposition \ref{general} to the case where $R=\tilde{\O}$, $k=4$, $\Delta$ is as above, $m$ is the number of roots of unity in $K$, and $\Omega$ be the set of all pairs $(a,b)=((a',a''),(b',b''))\in \tilde{\O}^2$ such that the following conditions hold true
\begin{enumerate}
	\item $a\tilde{\O}+b\tilde{\O}=\tilde{\O}$;
	
	\item Each of the elements $a'$,$a''$,$b'$,$b''$ is prime to $I$;
	
	\item $b'\O$ and $b''\O$ are prime ideals of $\O$  with residue characteristic prime to $m$;
	
	\item The residue characteristic of $b'\O$ does not divide $\eps_{\O}(b'')$, and vice versa.  
\end{enumerate}

Therefore it remains to verify conditions 1, 2 of Proposition \ref{general} in this setting. The following lemma verifies the condition 1.

{\lem \label{cond1} Under the conditions of Theorem \ref{main} any row $(a_1,b_1)$, where $a_1,b_1\in \tilde{\O}$ and $a_1\tilde{\O}+b_1\tilde{\O}=\tilde{\O}$, can be transformed to a row $(a^m,b)$, where $(a,b)\in\Omega$, by 4 elementary transvections from $\SL(2,\tilde{\O})$.}

\begin{proof}
	{\bf Step 1.} By one transvection, we transform the row $(a_1,b_1)$ into a row $(a_1,b_2)$, where $b_2$ is prime to the ideal $(mI,mI)\unlhd \tilde{\O}$. Equivalently, $b_2=(b_2',b_2'')$, where both $b_2'$ and $b_2''$ are prime to $mI$.
	
	Clearly, such a transformation is possible, because $\tilde{\O}/(mI,mI)$ is a semilocal ring.
	
	\medskip
	
	We shall make use of power norm residue symbols and the $m$-th power reciprocity law as described in \cite{BassMilnorSerre}, appendix on number theory. Let $p_1$,$\ldots$,$p_r$ be all the rational primes dividing $m$. Fix a positive integer $N$ that is so large that, firstly, for any $K$-prime $\P$ lying over one of the $p_i$ any integer element of the completion $K_{\P}$ that is congruent to 1 modulo $\P^N$ has an $m$-th root in $K_{\P}$; and secondly, powers in the decomposition of $I$ into $K$-primes do not exceed $N$. Further for any $p_i$ fix a $K$-prime $\P_i$ lying over $p_i$. Then we fix $\P_i$-local units $u_i$, $w_i$ such that  
	$$
	\zeta=\prod_{i=1}^r \rsymb{u_i,\,w_i}{\P_i}_m
	$$
	is a primitive $m$-th root of unity. Existence of such $u_i$ and $w_i$ follows from \cite{BassMilnorSerre} proposition A.17 and the standard properties of power norm residue symbols.
	\medskip
	
	{\bf Step 2.} By one transvection, we transform the row $(a_1,b_2)$ into a row $(a_2,b_2)=((a_2',a_2''), (b_2',b_2''))$, where
	
	\smallskip
	
	$\bullet$ $a_2'\equiv a_1'\mod b_2'\O$ and $a_2''\equiv a_1''\mod b_2''\O$;
	
	$\bullet$ $a_2'\O$ and $a_2''\O$ are prime ideals;
	
	$\bullet$ $a_2'\equiv a_2'' \equiv w_i\mod \P_i^N$ for $i=1$,$\ldots$,$r$;
	
	$\bullet$ $a_2'\equiv a_2'' \equiv 1\mod \frq^N$ for any prime $\frq$ such that $mI\sub \frq$, $\frq\ne\P_i$ $i=1$,$\ldots$,$r$;
	
	$\bullet$ $a_2''$ is positive in every real embedding of $K$;
	
	$\bullet$ If $K$ has a real embedding, then $a_2'$ is negative in one such embedding and positive in every remaining real embeddings.
	
	\smallskip
	
	Existence of such $a_2'$ and $a_2''$ follows from the Chinese remainder theorem and the
	generalized Dirichlet theorem on primes in arithmetic progressions, see A.11 in \cite{BassMilnorSerre}. By conditions on $a_2'$ and $a_2''$, we have $a_2'\equiv a_2''\mod I$, i.e. $a_2=(a_2',a_2'')\in \tilde{\O}$. Since $b_2'$ is prime to $I$, it follows by Lemma~\ref{quotient} that $a_2\equiv a_1 \mod b_2\tilde{\O}$.
	
	\medskip
	
	{\bf Step 3.} By one transvection, we transform the row $(a_2,b_2)$ into a row $(a_2,b)=((a_2',a_2''), (b',b''))$, where
	
	\smallskip
	
	$\bullet$ $b'\equiv b_2'\mod a_2'\O$ and $b''\equiv b_2''\mod a_2''\O$;
	
	$\bullet$ $b'\equiv b''\mod I$;
	
	$\bullet$ $b'\O$ and $b''\O$ are prime ideals;
	
	$\bullet$ $b'$ and $b''$ are prime to $mI$;
	
	$\bullet$ The residue characteristic of $b'\O$ does not divide $\eps_{\O}(b'')$, and vice versa;
	
	$\bullet$ $a_2'$ is an $m$-th power modulo $b'\O$, and $a_2''$ is an $m$-th power modulo $b''\O$.
	
	\smallskip
	
	It is enough to prove that such $b'$ and $b''$ exist; then it would follow by Lemma~\ref{quotient} that $b\equiv b_2\mod a_2\tilde{\O}$.
	
	It follows from the proof of Lemma~3 of \cite{CarterKeller} that there exists an element $v''\in \O$ prime to $m$ such that the congruence $b''\equiv v''\mod m^N\O$ guaranties that $a_2''$ is an $m$-th power modulo $b''\O$, provided $b''\O$ is a prime ideal and $b''\equiv b_2''\mod a_2''\O$. Since it only residue of $v''$ modulo $m^N$ that matters, we may assume that $v''\equiv 1 \mod \frq$ for any prime $\frq$ such that $I\sub \frq$ but $m\notin \frq$. Thus $v''$ is prime to $mI$. Also we may assume that $v''\equiv b_2''\mod a_2''\O$.
	
	Further let $\A=a_2''m^NI$ and $H_{\A}$ denote the ray class group of $K$ that correspond to the ideal $\A$. By the existence theorem of class field theory there is a finite abelian extension $K_{\A}/K$ for which the Artin reciprocity map gives an isomorphism $H_{\A}\simeq \Gal(K_{\A} /K)$.
	
	Now we choose $b'$ such that
	
	\smallskip
	
	$\bullet$ $b'\equiv b_2'\mod a_2'\O$;
	
	$\bullet$ $b'\equiv v''\mod I$;
	
	$\bullet$ $b'\O$ is a prime ideal;
	
	$\bullet$ The extension $K_{\A} /K$ is unramified in $b'\O$, i.e. $b'$ is prime to $\A$;
	
	$\bullet$ $b'$ is prime to $mI$;
	
	$\bullet$ $a_2'$ is an $m$-th power modulo $b'\O$;
	
	\smallskip 
	
	Let us prove that such $b'$ exists. Here we use the conditions of Theorem \ref{main}, and consider two cases.
	
	\smallskip 
	
	{\bf Case 1.} The field $K$ has no real embeddings and the ideal $I$ is prime to $m$.
	
	Similarly to $v''$ we can find $v'\in \O$ prime to $m$ such that the congruence $b'\equiv v'\mod m^N\O$ guaranties that $a_2'$ is an $m$-th power modulo $b'\O$, provided $b'\O$ is a prime ideal and $b'\equiv b_2'\mod a_2'\O$. Since $m\O$, $a_2'\O$ and $I$ are pairwise coprime there exist infinitely many prime elements $b'$ such that
	\begin{align*}
		&b'\equiv b_2'\mod a_2'\O\tc\\
		&b'\equiv v''\mod I\tc\\
		&b'\equiv v'\mod m^N\O\tp
		\end{align*}
	
	We chose one avoiding the primes that ramifies $K_{\A}$.
	
		\smallskip 
	
	{\bf Case 2.} The field $K$ has a real embedding.
	
	 In this case we have $m=2$. We search for a suitable $b'$ among prime elements such that
	\begin{align*}
		&b'\equiv b_2'\mod a_2'\O\tc\\
		&b'\equiv v''\mod m^NI\tp
	\end{align*}

	The only terms in the quadratic reciprocity law for $a_2'$ and $b'$ that have not been fixed yet are the one corresponding to the ideal $b'\O$ and the one corresponding to the real embedding that makes $a_2'$ negative. Therefore, we can guarantee that $a_2'$ is a square modulo $b'\O$ by prescribing certain sign for $b'$ in this embedding. That still leaves infinitely many prime elements. We chose one avoiding the primes that ramifies $K_{\A}$.
	
		\smallskip 
		
	Now let $p$ be the rational prime lying under $b'\O$. Since $K_{\A}/K$ is unramified in $b'O$ and $L_p$ is totally ramified in $b'O$ it follows that $K_{\A}\cap L_p=K$; hence $\Gal(L_pK_{\A}/K)=\Gal(L_p/K)\times \Gal(K_{\A}/K)$. We choose an element $\sigma\in \Gal(L_pK_{\A}/K)$ such that $\res_{K_{\A}}(\sigma)=(v''\O,K_{\A}/K)$, where $(v''\O,K_{\A}/K)$ is an Artin symbol, and $\res_{L_p}(\sigma)$ is nontrivial.
	
	By the Tchebotarev density theorem, there exists infinitely many primes $\P$ such that $(\P,L_pK_{\A}/K)=\sigma$. Choose one avoiding the ramified primes and the divisors of $\eps(b')mI$.
	
	Since $(\P,K_{\A}/K)=(v''\O,K_{\A}/K)$, it follows by Artin reciprocity that there exists $\lambda$ in $K$ such that it is multiplicatively congruent to $1 \mod\A$ and $\lambda v''\O = \P$. Set $b''=\lambda v''$. Then $b''\equiv v'' \mod\A$ and $b''\O=\P$. Since $(b''\O, L_p/K)$ is nontrivial, it follows that $\eps(b'')$ is not divisible by $p$, c.f. A.8 of \cite{BassMilnorSerre}.
	
	\medskip
	
	{\bf Step 4.} By one transvection, we transform the row $(a_2,b)$ into a row $(a^m,b)=((a',a'')^m, (b',b''))$, where $(a,b)\in\Omega$.
	
	By the previous step, there exist $a'$ and $a''$ such that $(a')^m\equiv a_2'\mod b'\O$ and $(a'')^m\equiv a_2''\mod b''\O$. Since $b'$ and $b''$ are prime to $I$, we may assume that $a'\equiv a''\equiv 1\mod I$; hence $(a,b)=((a',a''), (b',b''))\in\Omega$. Finally, by Lemma~\ref{quotient} we have $a^m\equiv a_2\mod b\tilde{\O}$.
 	\end{proof}

In order to verify the condition 2 of Proposition \ref{general}, we need the following lemma.

{\lem \label{cond2} In the setting of Theorem \ref{main}, let $\A$ be a nonzero ideal of $\O$, let $b\in \O$ be a nonzero element such that
\renewcommand{\theenumi}{\roman{enumi}}
\begin{enumerate}
\item  $b\O$ is a prime ideal with residue characteristic prime to $m$;

\item $b\O$ and $\A$ are comaximal. 
\end{enumerate}
\renewcommand{\theenumi}{\arabic{enumi}} 
Let $\S\sub \Pm$ be a finite set of rational primes such that it does not contain the residue characteristic of $b\O$ and one of the following conditions holds: either the Artin symbol $(b\O, L_p/K)$ is non trivial for all $p\in \S$; or the Artin symbol $(b\O, L_p/K)$ is trivial for all $p\in \S$, but there exists $\ovl{\sigma_1}\in \Gal(L/K)$, where $L$ is a composite $L=\prod_{p\in\S} L_p$, such that $\res_{L_p}(\sigma)$ is nontrivial for all $p\in \S$.

Then there exists $c\in \O$ such that $bc\equiv -1\mod \A$ and $\eps(c)=m\gamma$, where none of the primes from $\S$ divide $\gamma$. }
\begin{proof}
	Let $H_{\A}$ denote the ray class group of $K$ that correspond to the ideal $\A$. By the existence theorem of class field theory there is a finite abelian extension $K_{\A}/K$ for which the Artin reciprocity map gives an isomorphism $H_{\A}\simeq \Gal(K_{\A} /K)$. Consider the composite $L_{\A} = LK_{\A}$.
	
	Now we consider two cases described in the assumptions.
	
	{\bf Case 1} The Artin symbol $(b\O, L_p/K)$ is non trivial for all $p\in \S$.
	
	Applying the Tchebotarev density theorem to the extension $L_{\A} /K$, we conclude that there are infinitely many $K$-primes $\P$ such that $(\P,L_{\A}/K)=(b\O, L_{\A}/K)^{-1}$. Choose one avoiding the ramified primes and the divisors of $m\A$.
	
	Let $c_0$ be any element of $\O$ such that $bc_0\equiv -1 \mod \A$. Then $(b\O, K_{\A}/K)^{-1} = (c_0, K_{\A}/K)$ and hence $(\P, K_{\A}/K) = (c_0\O, K_{\A} /K)$. By Artin reciprocity it follows that there exists $\lambda$ in $K$ such that it is multiplicatively congruent to $1 \mod\A$ and $\lambda c_0\O = \P$. Set $c=\lambda c_0$. Then $bc\equiv -1 \mod\A$; and since $(c\O, L_{\A}/K)=(b\O, L_{\A}/K)^{-1}$ has nontrivial restrictions to all the $L_p$ for $p\in \S$, it follows that $\ord_p (\eps(c))=\ord_p(|R/\P|-1) = \ord_p(m)$ for all $p\in \S$. Therefore, $\eps(c)=m\gamma$, where none of the primes from $\S$ divide~$\gamma$.
	
	\medskip
	
	{\bf Case 2} The Artin symbol $(b\O, L_p/K)$ is trivial for all $p\in \S$, but there exists $\ovl{\sigma_1}\in \Gal(L/K)$ such that $\res_{L_p}(\sigma)$ is nontrivial for all $p\in \S$.
	
	Choose $\sigma_1\in \Gal(L_{\A}/K)$ such that $\res_{L}(\sigma_1)=\ovl{\sigma_1}$. Then set $\sigma_2=(b\O,L_{\A}/K)\sigma_1^{-1}$.
	
	Applying the Tchebotarev density theorem to the extension $L_{\A} /K$, we conclude that there are infinitely many $K$-primes $\P_1$, $\P_2$ such that $\sigma_i^{-1}=(\P_i, L_{\A}/K)$, for  $i = 1, 2$. Choose any two distinct such primes, avoiding the
	ramified primes and the divisors of $m\A$.
	
	Similarly to the previous case, we obtain that the ideal $\P_1\P_2$ is principal with the generator $c$ such that $bc\equiv -1\mod \A$. Since both $\sigma_1$ and $\sigma_2$ have nontrivial restrictions to all the $L_p$ for $p\in \S$, it follows that $\ord_p (|R/\P_1|-1)=\ord_p(|R/\P_2|-1) = \ord_p(m)$ for all $p\in \S$. Hence $\eps(c)=\lcm((|R/\P_1|-1),(|R/\P_2|-1))=m\gamma$, where none of the primes from $\S$ divide~$\gamma$.
	
	\end{proof}

Now we finish the proof of the main theorem. Recall that it remains to verify the condition 2 of proposition \ref{general} in the setting above.

So let $(a,b)=((a',a''),(b',b''))\in \Omega$. Since $b'$ is prime to $I$, it follows from Lemma~\ref{quotient} that
$$
\eps_{\tilde{\O}}(b)=\lcm(\eps_{\O}(b'),\eps_{\O}(b''))\tp
$$
Recall that $e_p=\ord_p(m)$. Set
$$
\S_0=\{p\in\Pm\colon \ord_p(\eps_{\tilde{\O}}(b)) > e_p\}\tp
$$

Note that $\S_0$ does not contain residue characteristics of $b'\O$ and $b''\O$. Indeed, all the primes from $\S_0$ are either divisors of $\eps(b')$ or of $\eps(b'')$. A divisor of $\eps(b')$ can not be the residue characteristic of $b''\O$ by the definition of $\Omega$ and it can not be the residue characteristic of $b'\O$, because $|\O/b'\O|=\eps(b')+1$. Similarly, a divisor of $\eps(b'')$ also can not coincide with these residue characteristics. Therefore, for $p\in\S_0$ the extension $L_p/K$ is unramified in $b'O$ and $b''\O$, and we can consider the Artin symbols $(b'\O,L_p/K)$ and $(b''\O,L_p/K)\in \Gal(L_p/K)$. 

By the definition of $\S_0$ for any $p\in \S_0$ we have either $\ord_p(\eps(b')) > e_p$ or $\ord_p(\eps(b'')) > e_p$. In the first case $(b'\O,L_p/K)$ is trivial, and in the second case $(b''\O,L_p/K)$ is trivial. Therefore, we have
$$
\S_0=\S_1\sqcup \S_2\sqcup \S_3\tc 
$$
where
\begin{gather*}
	\S_1=\{p\in\S_0\colon (b'\O,L_p/K)=1\text{ and } (b''\O,L_p/K)\ne 1 \}\tc\\
	\S_2=\{p\in\S_0\colon (b'\O,L_p/K)=(b''\O,L_p/K)=1 \}\tc\\
	\S_3=\{p\in\S_0\colon (b'\O,L_p/K)\ne 1\text{ and } (b''\O,L_p/K)=1 \}\tp
	\end{gather*}
 Now for $i=1$,$2$,$3$ set $\delta_i=|\S_i\cap \S_{\bad}|$, and $\Delta_i=\max(1,\ln(\delta_i+1)/\ln 2)$. Finally, set $\Delta_0=\Delta_1+\Delta_2+\Delta_3$. Clearly, $\Delta_0\le\Delta$.
 
 Let us prove that for any $i\in\{1,2,3\}$ we can present $\S_i$ as such a union
 $$
 \S_i=\bigcup_{j=1}^{\Delta_i} \S_i^{(j)}
 $$
 that for any $1\le j\le \Delta_i$ there exists an element $\sigma\in \Gal(L/K)$ with $\res_{L_p}(\sigma)$ being nontrivial for all $p\in \S_i^{(j)}$, where $L=\prod_{p\in\S_i} L_p$.
 
 Similarly to the proof of Lemma~4 of \cite{CarterKeller}, we obtain that
 $$
 \Gal(L/K)=\Gal(L_{\bad}/K)\times \prod_{p\in  \S_i \sm \S_{\bad}}\Gal(L_p/K) \tc
 $$
 where $L_{\bad}=\prod_{p\in \S_i \cap \S_{\bad}} L_p$.
 
 If $\delta_i=0$, then $\Delta_i=1$ and, clearly, we can set $\S_i^{(1)}=\S_i$. So assume that $\delta_i>0$. Then by Referee's Addendum to the paper \cite{CarterKeller}, there exist elements $\sigma_1$,$\ldots$,$\sigma_{\Delta_i}\in \Gal(L_{\bad}/K)$ such that for any $p\in\S_i\cap \S_{\bad}$ at least one of the $\sigma_j$ does not belong to the kernel of the restriction map $\Gal(\L_{\bad}/K)\to \Gal(L_p/K)$. Therefore we can set
 $$
 \S_i^{(j)}=\left(\S_i\sm \S_{\bad}\right)\cup \{p\in \S_i\cap\S_{\bad}\colon \res_{L_p}(\sigma_j)\ne 1\}\tp
 $$
 
 Now by Lemma~\ref{cond2} for any $1\le i\le 3$ and $1\le j\le \Delta_i$ there exist elements $c_{i,j}'$, $c_{i,j}''\in \O$ such that 
 \begin{align*}
 &b'c_{i,j}'\equiv -1\mod a'I\tc\\
  &b''c_{i,j}''\equiv -1\mod a''I\tc\\ 
  &\eps(c_{i,j}')=m\gamma_{i,j}'\tc\\
  &\eps(c_{i,j}''))=m\gamma_{i,j}''\tc
 \end{align*}
where none of the primes from $\S_i^{(j)}$ divides $\gamma_{i,j}'$ or $\gamma_{i,j}''$. 

Here we applied Lemma~\ref{cond2} for $\A=a'I$ resp. $\A=a''I$; and we used that $b'\O$ and $a'I$ are comaximal, because $b'\O$ and $a'\O$ are comaximal and $b'\O$ and $I$ are comaximal; similarly $b''\O$ and $a''I$ are comaximal.

Since the residue characteristics of $b'\O$ and $b''\O$ are prime to $m$, it follows that both $\eps(b')$ and $\eps(b'')$ are divisible by $m$. Let
$$\eps_{\tilde{\O}}(b)=\lcm(\eps(b'),\eps(b''))=m\gamma\tp$$

Now we set $c_{i,j}=(c_{i,j}',c_{i,j}'')$ and $\gamma_{i,j}=\gcd(\gamma,\lcm(\gamma_{i,j}',\gamma_{i,j}''))$. Let us verify that these elements satisfy the requirements in the condition 2 of Proposition \ref{general}.

\smallskip

$\bullet$ Since for all $i$,$j$ we have $c'_{i,j}\equiv c''_{i,j}\equiv -1\mod I$, it follows that $c_{i,j}\in \tilde{\O}$.

\smallskip

$\bullet$ All the primes dividing $\gamma_{i,j}$ also divides $\gamma$ and thus by definition of $\S_0$ they belong to $\S_0$. But the primes from $\S_i^{(j)}$ does not divide $\gamma_{i,j}'$ and $\gamma_{i,j}''$; hence they does not divide $\lcm(\gamma_{i,j}',\gamma_{i,j}'')$; hence they does not divide $\gamma_{i,j}$. Since $\S_0=\bigcup_{i,j} S_i^{(j)}$, it follows that $\gcd_{i,j}(\gamma_{i,j})=1$.

\smallskip

$\bullet$ Since $b'c_{i,j}'\equiv -1\mod a'\O$, $b''c_{i,j}''\equiv -1\mod a''\O$ and $a'$ is prime to $I$, it follows by Lemma~\ref{quotient} that $bc_{i,j}\equiv -1\mod a\tilde{\O}$.

\smallskip

$\bullet$ Finally, since $c_{i,j}'$ is prime to $I$, it follows by Lemma~\ref{quotient} that $\eps_{\tilde{\O}}(c_{i,j})=\lcm(\eps_{\O}(c_{i,j}',\eps_{\O}(c_{i,j}'')))$. Therefore,
$$
\gcd(\eps_{\tilde{\O}}(b),\eps_{\tilde{\O}}(c_{i,j}))=\gcd(m\gamma,\lcm(m\gamma_{i,j}',m\gamma_{i,j}''))=m\gcd(\gamma,\lcm(\gamma_{i,j}',\gamma_{i,j}''))=m\gamma_{i,j}\tp
$$
That finishes the proof of Theorem \ref{main}.

\section{Final remarks}
\label{finalremarks}

Now we give a quick proof that $E(n,\O_S,I)$ has finite width in Tits--Vaserstein generators without any assumptions on $K$ and $I$ that relies only on the result of \cite{CarterKeller}. This proof is similar to the one given in \cite{TavgenThesis} for $E(\Phi,I)$; and we will see why this proof does not allow to obtain any explicit estimate.

{\prop In the previous notation, let $I$ be a non zero ideal in $\O_S$, and let $n\ge 3$. Then the group $E(n,\O_S,I)$ has finite width in Tits--Vaserstein generators.}
\begin{proof}
	For $g\in E(n,\O_S,I)$ let $l(g)$ be the smallest integer such that $g$ can be expressed as a product of $l(g)$ Tits--Vaserstein generators.
	
	Note that $\O_S/I$ is a finite set, so let $\zeta_1$,$\ldots$,$\zeta_r$ be representatives of all the classes in $\O_S/I$. By \cite{CarterKeller} any element $g\in E(n,\O_S,I)\le E(n,\O_S)$ is a product of at most $N$ transvections. We can write it as
	$$
	g=\prod_{u=1}^N t_{i_u,j_u}(\zeta_{k_u}+\xi_u)\tc
	$$  
	where $\xi_u\in I$. Then we have
	$$
	g=\prod_{u=1}^N t_{i_u,j_u}(\xi_u)^{h_u}\cdot \prod_{u=1}^N t_{i_u,j_u}(\zeta_{k_u})\tc
	$$
	where $h_u=\prod_{v=1}^u t_{i_u,j_u}(-\zeta_{k_v})$. The number $l(t_{i_u,j_u}(\xi_u)^{h_u})$ can be estimated in a manner similar to the proof of Corollary~\ref{TitsVaserstein} with some constant $C$. Now let $\Theta$ be the set of all elements of the form
	$$
	\prod_{u=1}^N t_{i_u,j_u}(\zeta_{k_u})\tc
	$$
	then $\Theta$ is a finite set. Therefore we have.
	$$
	l(g)\le CN+\max_{h\in\Theta\cap E(n,\O_S,I)} l(h)\tp
	$$
\end{proof}
As we can see this estimate is not explicit because we do not know the bound for $l(h)$, where $h\in \Theta\cap E(n,\O_S,I)$. We just know that some way to express such $h$ as a product of Tits--Vaserstein generators exists.

However, we believe that it is possible to generalise the results of the present paper to the arbitrary ideals in $\O_S$, replacing $\SL(n,\O_S,I)$ by $E(n,\O_S,I)$. The only problem is the compatibility of certain congruences in step 3 of the proof of Lemma~\ref{cond1}. Our conjecture is that those congruences will be compatible automatically if the matrix we start with belong to $E(3,\tilde{\O_S})$.  

\end{document}